\documentclass[12pt]{article}
\usepackage{graphicx}
\usepackage{amsmath,amsthm,amssymb,enumerate}
\usepackage{euscript,mathrsfs}
\usepackage{color}
\usepackage{dsfont}
\usepackage{url}
\usepackage[notref,notcite]{showkeys}
\usepackage[left=2cm,right=2cm,top=3.5cm,bottom=3.5cm]{geometry}
\usepackage{color}
\usepackage[framemethod=tikz]{mdframed}
\allowdisplaybreaks

\usepackage{soul}

\catcode`\@=11 \@addtoreset{equation}{section}

\catcode`\@=12

\newtheorem{Theorem}{Theorem}[section]
\newtheorem{Proposition}[Theorem]{Proposition}
\newtheorem{Lemma}[Theorem]{Lemma}
\newtheorem{Corollary}[Theorem]{Corollary}

\theoremstyle{definition}
\newtheorem{Definition}[Theorem]{Definition}

\newtheorem{Remark}[Theorem]{Remark}

\newcommand{\bTheorem}[1]{
	\begin{Theorem} \label{T#1} }
	\newcommand{\eT}{\end{Theorem}}

\newcommand{\bProposition}[1]{
	\begin{Proposition} \label{P#1}}
	\newcommand{\eP}{\end{Proposition}}

\newcommand{\bLemma}[1]{
	\begin{Lemma} \label{L#1} }
	\newcommand{\eL}{\end{Lemma}}

\newcommand{\bCorollary}[1]{
	\begin{Corollary} \label{C#1} }
	\newcommand{\eC}{\end{Corollary}}

\newcommand{\bRemark}[1]{
	\begin{Remark} \label{R#1} }
	\newcommand{\eR}{\end{Remark}}

\newcommand{\bDefinition}[1]{
	\begin{Definition} \label{D#1} }
	\newcommand{\eD}{\end{Definition}}

\newcommand{\Ds}{\mathbb{D}_x}

\newcommand{\bu}{\mathbf u}

\newcommand{\bfphi}{\boldsymbol{\varphi}}

\newcommand{\bFormula}[1]{
	\begin{equation} \label{#1}}
	\newcommand{\eF}{\end{equation}}

\newcommand{\vrn}{\vr_n}
\newcommand{\vun}{\vu_n}

\newcommand{\Ov}[1]{\overline{#1}}

\newcommand{\aleq}{\stackrel{<}{\sim}}

\newcommand{\vr}{\varrho}

\newcommand{\vu}{\vc{u}}

\newcommand{\vc}[1]{{\bf #1}}

\newcommand{\Div}{{\rm div}_x}
\newcommand{\Grad}{\nabla_x}

\newcommand{\dx}{\,{\rm d} {x}}

\newcommand{\dt}{\,{\rm d} t }

\newcommand{\intO}[1]{\int_{\Omega} #1 \ \dx}

\newcommand{\D}{{\rm d}}

\newcommand{\ep}{\varepsilon}

\newcommand{\br}{ \nonumber \\ }

\def\softd{{\leavevmode\setbox1=\hbox{d}%
		\hbox to 1.05\wd1{d\kern-0.4ex{\char039}\hss}}}
\definecolor{Cgrey}{rgb}{0.85,0.85,0.85}
\definecolor{Cblue}{rgb}{0.50,0.85,0.85}
\definecolor{Cred}{rgb}{1,0,0}
\definecolor{fancy}{rgb}{0.10,0.85,0.10}
\definecolor{amaranth}{rgb}{0.9, 0.17, 0.31}

\newcommand\Cbox[2]{%
	\newbox\contentbox%
	\newbox\bkgdbox%
	\setbox\contentbox\hbox to \hsize{%
		\vtop{
			\kern\columnsep
			\hbox to \hsize{%
				\kern\columnsep%
				\advance\hsize by -2\columnsep%
				\setlength{\textwidth}{\hsize}%
				\vbox{
					\parskip=\baselineskip
					\parindent=0bp
					#2
				}%
				\kern\columnsep%
			}%
			\kern\columnsep%
		}%
	}%
	\setbox\bkgdbox\vbox{
		\color{#1}
		\hrule width  \wd\contentbox %
		height \ht\contentbox %
		depth  \dp\contentbox
		\color{black}
	}%
	\wd\bkgdbox=0bp%
	\vbox{\hbox to \hsize{\box\bkgdbox\box\contentbox}}%
	\vskip\baselineskip%
}

\mdfdefinestyle{MyFrame}{%
	linecolor=black,
	outerlinewidth=1pt,
	roundcorner=5pt,
	innertopmargin=\baselineskip,
	innerbottommargin=\baselineskip,
	innerrightmargin=10pt,
	innerleftmargin=10pt,
	backgroundcolor=white!20!white}

\newcommand{\mc}{\mathcal}

\newcommand{\T}{\mathbb{T}}

\allowdisplaybreaks

\begin{document}


\title{\bf On the long time behaviour of a system of several rigid bodies immersed in a viscous fluid}

\author{Marco Bravin \thanks{The work of M.B. was supported by the project "Ecuaciones en derivadas parciales motivadas por procesos de difusión y mecánica de fluidos: propiedades, asintóticas y homogeneización" (Ayuda financiada contrato Programa Gob. Cantabria -UC). The work of E.F. was partially supported by the		Czech Sciences Foundation (GA\v CR), Grant Agreement		24--11034S. The Institute of Mathematics of the Academy of Sciences of		the Czech Republic is supported by RVO:67985840. 	E.F. is a member of the Ne\v cas Center for Mathematical Modelling. A.R is supported by the Grant RYC2022-036183-I funded by MICIU/AEI/10.13039/501100011033 and by ESF+. A.R, A.Z have been partially supported by the Basque Government through the BERC 2022-2025 program and by the Spanish State Research Agency through BCAM Severo Ochoa CEX2021-001142-S and through project PID2023-146764NB-I00 funded by MICIU/AEI/10.13039/501100011033 and cofunded by the European Union. A.Z. was also partially supported  by a grant of the Ministry of Research, Innovation and Digitization, CNCS - UEFISCDI, project number PN-III-P4-PCE-2021-0921, within PNCDI III.} 
\and Eduard Feireisl
 \and Arnab Roy$^{1,2}$ \and Arghir Zarnescu$^{1,2,3}$
}

\date{}

\maketitle

\centerline{Universidad de Cantabria, E.T.S. de Ingenieros Industriales y de Telecomunicacíon}

\centerline{Departmento de Matemática Aplicada y Ciencias de la Computacíon}

\centerline{Avd. Los Castros 44, 39005 Santander, Spain}

\bigskip


\centerline{\v Zitn\' a 25, CZ-115 67 Praha 1, Czech Republic}

\centerline{$^1$ BCAM, Basque Center for Applied Mathematics}

\centerline{Mazarredo 14, E48009 Bilbao, Bizkaia, Spain.}

\centerline{$^2$IKERBASQUE, Basque Foundation for Science, }

\centerline{Plaza Euskadi 5, 48009 Bilbao, Bizkaia, Spain.}

\centerline{$^3$``Simion Stoilow" Institute of the Romanian Academy,}

\centerline{21 Calea Grivi\c{t}ei, 010702 Bucharest, Romania.}

\begin{abstract}
	
We consider several rigid bodies immersed in a viscous Newtonian fluid 
contained in a bounded domain in $R^3$. We introduce a new concept 
of \emph{dissipative} weak solution of the problem based 
on a combination of the approach proposed by Judakov with a suitable form of energy inequality. We show that global--in--time dissipative solutions always exist as long as the rigid bodies are connected compact sets. In addition, in the absence of external driving forces, the system always tends to a static equilibrium as time goes to infinity. The results hold independently of possible collisions of rigid bodies and for any finite energy initial data.

\end{abstract}


{\small


\medbreak


}

\section{Introduction}
\label{i}

Recently, the authors in \cite{ErvMaiTuc, MaiTuc} studied the long time behaviour of a single rigid body immersed in a viscous fluid filling the entire space $R^3$. Specifically, if the initial velocity of the system is small enough, 
the body approaches some final configuration as time goes to infinity. 
As pointed out in \cite{MaiTuc}, the methodology could be adapted to the system of several rigid bodies provided the initial velocities are so small that their mutual contacts are excluded. Their approach is based on the concept of strong solution to the fluid structure system, therefore 
small data are necessary in view of the well known regularity issues related to 
the 3D incompressible Navier--Stokes system. Let us mention some available works related to the long time behaviour of the fluid-rigid body system in the literature. The authors in \cite{BJ} analyzed the asymptotic behavior of a rigid ball interacting with an incompressible viscous fluid, establishing precise decay profiles for the coupled system in $R^3$. In \cite{EHL}, the long-time dynamics of a rigid disk moving in a two-dimensional viscous fluid, proving convergence toward steady states under small initial data has been studied. In \cite{BHPS}, the authors investigated the interaction of an anisotropic rigid body with a Poiseuille flow in an unbounded two-dimensional channel, characterizing the asymptotic velocity and orientation of the body over time. In \cite{FH}, the authors constructed solutions with unbounded energy for the fluid–disk system and analyzed the long-time behavior for large initial data, highlighting nonlinear phenomena absent in small-data regimes. The article \cite{GP} is devoted to the analysis of large time behavior of a rigid body of arbitrary shape under the presence of external forces and torques and it shows that if the external forces vanish for large time in $L^2$ sense then the global strong solution converges to zero as time approaches infinity. The large time behavior of solutions around nontrivial basic states for fluid-rigid ball interaction problem has been analyzed in \cite{TH}. In \cite{GPPM}, the authors establish that if the initial data for the fluid-rigid body system is small under appropriate norm, then the kinetic energy of the system vanishes as $t\rightarrow \infty$. Related results on the long-time behaviour of fluid-rigid body interaction systems under a spring-damper feedback force were obtained in \cite{takahashi2015stabilization, maity2023global, roy2021stabilization}.

We propose an alternative approach  based on an appropriate concept of weak solutions. For the system of equations that describes the evolution of a family of rigid bodies immersed in a Navier-Stokes fluid contained in a bounded domain $\Omega \subset R^3$, we introduce a concept of \emph{dissipative} weak solution based on the original idea of Judakov \cite{Juda} combined with a suitable form of energy balance. These weak solutions  naturally include the case of strong solutions while having some more structure than the distributional solutions, namely having what can be described informally as a differential form of the energy inequality. We show 
that global--in--time dissipative solutions always exist for any finite energy initial data. Moreover, in the absence of driving forces, 
the system stabilizes to a static equilibrium for large time. 
The result holds for 
\begin{itemize}
	\item arbitrarily large initial kinetic energy of the system (no smallness assumption on the velocity); 
	\item independently of possible collisions of several rigid bodies and/or collisions of the bodies with the domain boundary.
\end{itemize}

The above scenario should be compared with that of the Navier-Stokes fluid 
without rigid bodies and a possibly large initial velocity. Although the solution may not be regular, the velocity stabilizes exponentially to zero as long as external forces are absent.

\subsection{Problem formulation}
{
We consider a {domain} $\Omega \subset R^3$ containing a viscous incompressible fluid and with several rigid bodies  moving inside the fluid. The rigid bodies are connected compact sets $(\mathcal{S}^{i})_{i = 1}^N$. The {part of the domain $\Omega$ occupied by the fluid} at a time $t$ is given by
\[
\Omega_{F}(t) = \Omega \setminus \cup_{i = 1}^N \mc{S}^i(t),\ \mbox{and {we} set}\ Q_{F} = 
\left\{ (t,x) \ \Big|\ t \in (0,T),\ x \in \Omega_{F}(t) \right\}.
\] For the sake of simplicity, we set the fluid density $\vr_F$ to be $1$. {The motion of the fluid in $Q_{F}$} is described by
the fluid velocity $\vu_F = \vu_{F}(t,x)$ satisfying the \emph{Incompressible Navier--Stokes system}:
\begin{align} 
	\Div \bu_{F} &= 0, \label{m1} \\ 	
	\partial_t \bu_{F} + \Div (\bu_{F} \otimes \bu_{F}) + \Grad \Pi &= {\Div \mathbb{S} (\Ds \bu_{F})} + \vc{g}, \label{m2} \\ 
	{\mathbb{S}(\Ds \bu) = \mu \Ds \bu },\ 
	\Ds \bu& = \frac{ \Grad \bu + \Grad^t \bu }{2} ,\ \mu > 0 \label{m3}, 
\end{align}
where $\Pi$ is the pressure and the function $\vc{g}$ denotes an external {volume} force. Denoting $\vr_{\mathcal{S}^i}$ the mass density of each body, we define the total mass $m^{i}$ and the center of mass $\vc{h}_{i}(t)$ as
\begin{equation*}
m^i = \int_{\mc{S}^i(t)} \vr_{\mathcal{S}^i}(t,x)\, {\rm d} x,\quad \vc{h}_{i}(t)= \frac{1}{m^i} \int_{\mc{S}^i(t)} \vr_{\mathcal{S}^i}(t,x)\, x\, {\rm d} x.
\end{equation*}
Moreover, we consider the inertial tensor $J^{i}$ defined through its action on vectors $a$ and $b$ as:
 \begin{equation*}
J^i a \cdot b
= \int_{\mc{S}^i(t)} \vr_{\mathcal{S}^i}(t,x)\, [a \times (x - \vc{h}_{i}(t))] \cdot [b \times (x - \vc{h}_{i}(t))]\, {\rm d} x.
\end{equation*}
The associated ``rigid'' velocities at each point $x\in \mc{S}^{i}(t)$ are given by 
\begin{equation} \label{m5}	 
\bu_{\mc{S}^i} (t,x) = \vc{Y}_i(t)  + \mathbb{Q}_{i}(t) ( x - \vc{h}_i(t)),
\end{equation}
where
\[
	\vc{Y}_i(t)  = \frac{\D }{\dt } \vc{h}_i (t) ,\ \mathbb{Q}_i(t) = 
	\frac{ \D }{\dt} \mathbb{O}_i (t) \circ {\mathbb{O}_i}^{-1}(t),\ \mathbb{O}_i \in SO(3).
\]
Here $\vc{Y}_i$ denotes the translational velocity and $\mathbb{Q}_{i}$ is the angular velocity of the body. We suppose $\vc{Y}_i \in L^\infty_{\rm loc}([0, \infty); R^3)$, $\mathbb{Q}_i \in L^\infty_{\rm loc} ([0, \infty); R^{3 \times 3})$, 
$i=1,\dots, N$. The matrix $\mathbb{Q}_{i}$ is skew-symmetric, therefore it can be represented by a vector $\omega^i$, 
\begin{equation*}
\mathbb{Q}_{i}(t)\big(x - \vc{h}_i(t)\big) = \omega^i \times (x - \vc{h}_i(t)).
\end{equation*}
The balance of linear and angular momentum for the body $\mc{S}^i(t)$ reads
\begin{align}
m^i \frac{d}{dt} \vc{Y}_i(t)
&= \int_{\partial \mc{S}^i(t)} (\mathbb{S} (\Ds \bu_{F}) - pI)\, n\, d\sigma
+ \int_{\mc{S}^i(t)} \vr_{\mathcal{S}^i} \vc{g}\, {\rm d} x, \label{rigid1}\\[1ex]
J^i(t) \frac{d}{dt} \omega^i(t)
&= J^i \omega^i(t) \times \omega^i(t)
+ \int_{\partial \mc{S}^i(t)} (x - \vc{h}_i(t)) \times (\mathbb{S} (\Ds \bu_{F}) - pI)\, n\, d\sigma \notag \\
&\quad + \int_{\mc{S}^i(t)} \vr_{\mathcal{S}^i} (x - \vc{h}_i(t)) \times \vc{g}\, {\rm d} x,
\qquad i=1,\dots,N. \label{rigid2}
\end{align}
The configuration of the rigid bodies at a time $t \geq 0$ is parametrized by a family of affine isometries $(\sigma_i(t))_{i=1}^N$, such that 
\begin{equation} \label{m4}
\mc{S}^i (t) = \sigma^i [\mc{S}^i],\ \sigma^i(t)x = 
\mathbb{O}_{i} (t) x + \vc{h}_{i} (t),\ \mathbb{O}_i \in SO(3).
\end{equation}



The state of the system is described by its mass density 
$\vr = \vr(t,x)$ and the velocity $\vu = \vu(t,x)$. These quantities are defined throughout $\Omega$ by prescribing their values in the fluid region and inside each rigid body, namely:
\begin{align} 
\vr(t, x) = \left\{ \begin{array}{l} 1 (= \vr_F) \ \mbox{for}\ 
	x \in \Omega_F(t), \\ 
	\vr_{\mathcal{S}^i} \ \mbox{for}\ x \in \mathcal{S}^i(t),
	\end{array} \right. 
	\label{m6} \\
\vu(t,x) = 	\left\{ \begin{array}{l} \vu_F(t,x) \ \mbox{for}\ 
	x \in \Omega_F(t), \\ 
	\vu_{\mathcal{S}^i}(t,x) \ \mbox{for}\ x \in \mathcal{S}^i(t) 
\end{array} \right. \label{m7}
\end{align}
for a.a. $t > 0$. Moreover we assume {no-slip boundary condition}:

\begin{equation}
{\vu_F |_{\partial \mathcal{S}_i} = \vu_{\mathcal{S}_i},\ i = 1,\dots, N} ,\ \vu|_{\partial \Omega} = 0. 
\end{equation}
For the sake of simplicity, we have supposed the mass densities 
$\vr_{\mathcal{S}^i}$ of the rigid bodies are constant.
}
\medskip

\subsubsection{Mass conservation} 

The density $\vr$ and the velocity $\vu$ satisfy the equation of continuity 
\begin{equation} \label{m8}
\int_0^\infty \intO{ \Big[ 
	\vr \partial_t \varphi + \vr \vu \cdot \Grad \varphi \Big] } \dt
	= - \intO{ \vr_0 \varphi (0, \cdot) }	
\end{equation}	
for any $\varphi \in C^1_c([0, \infty) \times \Ov{\Omega})$.
As the velocity $\vu$ is solenoidal, DiPerna-Lions theory \cite{DL} 
{asserts that for a given velocity field} 
\begin{equation} \label{m10}
		\vu \in L^\infty ([0, \infty); L^2(\Omega; R^3)) \cap 
		L^2(0, \infty; W^{1,2}_0 (\Omega; R^3)),\ \Div \vu = 0,
\end{equation}
{equation \eqref{m8} admits 
a renormalized solution $\vr$ unique in the class} 
\begin{equation} \label{m9}
\vr \in C([0, \infty; L^1(\Omega)) \cap L^\infty ((0, \infty) \times \Omega). 
\end{equation}

\subsubsection{Momentum balance}

The weak formulation of the momentum balance goes back to Judakov \cite{Juda} with a slight modification introduced in \cite{EF64, EF65}. 
We require the integral identity
\begin{align} 
\int_0^\infty &\intO{ \Big[ \vr \vu \cdot \partial_t \bfphi + \vr \vu \otimes \vu : \Grad \bfphi - \mathbb{S}(\Ds \vu) : \Ds \bfphi  \Big] } \dt \br &= - \int_0^\infty \intO{ \vr \vc{g} \cdot \bfphi } \dt- 
\intO{ (\vr_0 \vu_0) \cdot \bfphi(0, \cdot)}
\label{m11}
\end{align}
to hold for any test function $\bfphi \in C^1_c([0, \infty) \times \Omega; R^3)$, $\Div \bfphi = 0$,
\begin{equation} \label{m12}
\Ds \bfphi (t, \cdot) = 0 \ \mbox{on an open neighbourhood of}\ 
\mathcal{S}^i(t),\ i = 1,\dots, N, 
\end{equation}
{
where the velocity field $\vu$ belong to the class \eqref{m10} and satisfies the compatibility conditions 
\eqref{m7}. In particular,}
\[
\Ds \bfphi (t, \cdot) = 0 \ \mbox{on}\ \mathcal{S}^i (t) \ \mbox{for a.a.}\ t \in (0, \infty).
\]
{Note that the weak formulations \eqref{m8}, \eqref{m11}, together with the compatibility condition \eqref{m7}
includes both the Navier-Stokes system satisfied in the fluid region as well as the validity of the momentum equations \eqref{rigid1}, 
\eqref{rigid2} on the rigid bodies.}

\subsubsection{Energy balance}

In contrast with \cite{EF64,EF65}, 
we require the energy balance to hold in the differential form
\begin{align} 
- \int_0^\infty \partial_t \psi \intO{ \frac{1}{2} \vr |\vu|^2 } \dt
&+ \int_0^\infty \psi \intO{ \mathbb{S}(\Ds \vu) : \Ds \vu } \dt \br  
&\leq \frac{1}{2} \psi (0) \intO{ \vr_0 |\vu_0|^2 }
+ \int_0^\infty \psi \intO{ \vr \vc{g} \cdot \vu } \dt
\label{m13}
\end{align} 	
for any $\psi \in C^1_c[0, \infty)$, $\psi \geq 0$.

\subsection{Dissipative weak solutions}
\label{dws}

\begin{Definition}[\bf Dissipative weak solution] \label{D1}

We say that the quantities 
\[ 
(\vr, \vu) \ \mbox{and}\ (\mathcal{S}^1(t, \cdot), \dots, 
\mathcal{S}^N (t, \cdot) )_{t \geq 0}
\]
represent \emph{dissipative solution} of the fluid--structure interaction problem if 
\begin{itemize}
	\item
the motion of the rigid bodies is determined by \eqref{m4}; 
\item the density $\vr$ and the velocity $\vu$ {belong to the class \eqref{m10}, \eqref{m9}} and satisfy the compatibility conditions \eqref{m6}, \eqref{m7}; 
\item the integral equalities \eqref{m8}, \eqref{m11} as well as the 
energy inequality \eqref{m13} hold.
\end{itemize}	
\end{Definition}	

\noindent
The original concept of weak solution \` a la Judakov is based on a more restrictive interpretation of the momentum equation, where the test functions satisfy the 
compatibility condition $\Ds \bfphi = 0$ only on the rigid bodies. 
Enlarging the class of test functions to \eqref{m12} allows to accommodate possible contacts of several rigid bodies and/or rigid bodies with the boundary $\partial \Omega$. Such a concept of weak solution was used in \cite{EF64}, \cite{EF65} to develop the theory of global in time weak solutions including possible collisions in both 
the compressible and the incompressible case. In both \cite{EF64} 
and \cite{EF65}, the energy 
inequality was imposed in the following form 
\begin{align} 
	\intO{ \frac{1}{2} \vr |\vu|^2 (\tau, \cdot) } +
	& \int_0^\tau \intO{ \mathbb{S}(\Ds \vu) : \Ds \vu } \dt \br  
	&\leq \intO{ \frac{1}{2} \vr_0 |\vu_0|^2 }
	+ \int_0^\tau \intO{ \vr \vc{g} \cdot \vu } \dt
	\label{m14}
\end{align}
for a.a. $\tau \geq 0$.

Our goal is to show that the construction of weak solutions presented in \cite{EF65} can be modified to accommodate the differential 
version of the energy inequality \eqref{m13}. This has an immediate impact on the long--time behaviour of global in time solutions. Whether or not the same is possible for the rigid bodies immersed in a viscous \emph{compressible} fluid remains an open problem. 

The paper is organized as follows. In Section \ref{r}, we revisit 
the construction of the weak solution in \cite{EF65} and show how to 
accommodate the differential form of the energy inequality \eqref{m13}. 
Applications to the long--time behaviour of dissipative solutions are discussed in Section \ref{L}.

\section{Weak solutions to the fluid structure interaction problem revisited}
\label{r} 

Our goal is to show the existence of global--in--time dissipative solutions for any finite energy initial data $\vr_0$, $\vu_0$ and 
any distribution of the rigid objects 
\begin{equation} \label{rr1}
\mathcal{S}^i,\ \mathcal{S}^i = {\rm cl}[{\rm int}[\mathcal{S}^i]],\ 
\mathcal{S}^i \ \mbox{compact, connected},\ i = 1,\dots, N
\end{equation}
such that their initial positions satisfy 
\begin{equation} \label{rr2}
{\rm int}[\mathcal{S}^i(0)] \cap {\rm int}[\mathcal{S}^j(0)] = \emptyset 
\ \mbox{as soon as}\  i \ne j.
	\end{equation}
To this end, 
we follow step by step the construction presented in 
\cite{EF65} and show how to accommodate the strong form of the energy inequality \eqref{m13}. 

The strategy of \cite{EF65} consists of several steps:
\begin{enumerate} 
	\item
For smooth and initially well separated bodies and a smooth spatial domain $\Omega$ construct a weak solution up to the time of the first contact, meaning either the contact of two or more rigid bodies or a contact of 
a rigid body with $\partial \Omega$, see Proposition 2.3 and 
Corollary 2.4 in \cite{EF65}.

\item Extend the result to arbitrary compact bodies as in 
\eqref{rr1} and $\Omega$ a general bounded open set, see   \cite[Proposition 3.1]{EF65}. This is done by means of an approximation procedure using the result of Step 1. The solution exist up to the time of the first contact. 

\item Continue the solution after the time of the first contact considering the colliding rigid objects as a single body and repeating 
the arguments of Step 2.
As there is a finite number of possible collision times, the process gives rise to a global in time weak solution. This has been explored in \cite[Section 4]{EF65}.

\end{enumerate}	

It is worth noting that at each step of the above procedure, one has to perform the limit in the convective term 
\[
\vrn \vun \otimes \vun \to \vr \vu \otimes \vu 
\]
in a weak sense. In particular, one obtains the convergence of the
weak time derivative of the kinetic energy 
\[ 
\int_0^\infty \partial_t \psi \intO{ \frac{1}{2} \vr_n |\vu_n|^2 }\dt 
\to \int_0^\infty \partial_t \psi \intO{ \frac{1}{2} \vr |\vu|^2 }\dt. 
\]
Thus the differential form of the energy inequality \eqref{m13} remains 
valid for the limit solution as long as it holds for the approximate solutions 
constructed in Step 1. Our goal therefore reduces to finding a suitable 
construction of solutions with smooth data (up to the first contact) that would satisfy the energy balance \eqref{m13}. 

\subsection{Approximate solutions with smooth data}

In accordance with the previous discussion, we consider the situation 
of Step 1 of initially well separated smooth bodies and a smooth spatial domain $\Omega$. Our goal is construct a weak solution up to the first contact satisfying the differential form of the energy inequality 
\eqref{m13}.

We use the idea of San Martin, Starovoitov, and Tucsnak \cite{SST} and approximate the rigid bodies by a fluid of high viscosity. 
Specifically, adapting the strategy proposed in \cite{FeHiNe} we proceed as follows:
\begin{itemize}
\item 
As $\Omega \subset R^3$ is bounded, we may suppose 
\[
\Omega \subset (-L,L)^3 
\]
for $L > 0$ large enough. Thus we may consider the problem with the 
spatially periodic conditions, meaning on the flat torus 
\[
\T^3 = \left( [-L,L]|_{\{ -L; L \} } \right)^3.
\] 

\item We consider a system of equations 

{
\begin{align} 
\partial_t \vr + \Div (\vr [\vu]_\delta ) &= 0,	
	\label{r1} \\
\partial_t (\vr \vu) + \Div (\vr \vu \otimes [\vu]_\delta) + \Grad P &= 
\Div ([\mu_{\ep}]_\delta \Ds \vu ) - \chi_\ep \vu + \vr \vc{g}, 
\label{r2} \\
\partial_t \mu_{\ep} + \Div (\mu_{\ep} [\vu]_\delta) &= 0,\ \label{r3}\\ 
\Div \vu &= 0, \label{r4}
\end{align}
}
where 
\[
[h]_\delta = \sigma_\delta * h 
\]
denotes a convolution in space with a family of regularization kernels with 
{
\begin{equation*}
\sigma_{\delta}(x)=\frac{1}{\delta^3}\sigma\left(\frac{|x|}{\delta}\right),\ \sigma\in C_c^{\infty}(-1,1),\ \sigma(z)>0 \mbox{ for }-1<z<1, \sigma(z)=\sigma(-z),\ \int\limits_{-1}^{1} \sigma(z)\ dz=1.
\end{equation*}
}

The quantity $\chi_\ep$ represents a  
penalization term, 
\[
\chi_\ep = \frac{1}{\ep} \chi,\ \chi \in C^\infty(\T^3),\ 
\chi(x) = \left\{ \begin{array}{l} 0 \ \mbox{for}\ x \in \Ov{\Omega}, \\ 
> 0 \ \mbox{for}\ x \in \T^3 \setminus \Ov{\Omega}. \end{array} \right. 	
\]

{The quantity $\mu_\varepsilon$ provides a regularized characterization of the rigid body regions, whose evolution is governed by the transport equation \eqref{r3}; {similarly the  friction coefficient $\chi_\ep$ is a penalization of the no-slip boundary condition on 
	$\partial \Omega$},	 see \cite{FeHiNe} for details.
}

\item System \eqref{r1}--\eqref{r4} is supplemented with smooth initial data
\begin{align}
\vr(0, \cdot) = \vr_{0, \delta} = 
\left\{ \begin{array}{l} 1 (= \vr_F) \ \mbox{for}\ x \in \T^3 
\setminus \cup_{i=1}^N \mathcal{S}^i (0)	,\\ 
\in [0, \vr_{S^i}] \ \mbox{for}\ x \in \mathcal{S}^i(0),\ 
{\rm dist}[x, \partial \mathcal{S}^i(0)] < \delta, \\ 
\vr_{S^i} \ \mbox{for}\ x \in \mathcal{S}^i(0),\ 
{\rm dist}[x, \partial \mathcal{S}^i] \geq \delta,\ i = 1,\dots, N;
\end{array}\right. 		
	\label{r5}	
	\end{align}
	\begin{align}
		\vu(0, \cdot) = \vu_{0, \delta} = 
		\left\{ \begin{array}{l} 0 \ \mbox{for}\ x \in \T^3 
			\setminus \Omega	,\\
			\vu_{0,F} \ \mbox{for}\ x \in \Omega \setminus \cup_{i=1}^N 
			\mathcal{S}^i(0),\\ 
			\vu_{0,S^i} \ \mbox{for}\ x \in \mathcal{S}^i(0);
		\end{array}\right. 		
		\label{r6}	
	\end{align}
	\begin{align}
	\mu(0, \cdot) = \mu_{0, \ep, \delta} = 
	\left\{ \begin{array}{l} 1 \ \mbox{for}\ x \in \T^3 
		\setminus \cup_{i = 1}^N \mathcal{S}^i	,\\
		\in \left[ 0,\ \frac{1}{\ep} \right] \ \mbox{for}\ 
		x \in \mathcal{S}^i(0),\ {\rm dist}[x, \partial \mathcal{S}^i] < \delta, \\
		\frac{1}{\ep} \ \mbox{for}\ x \in \mathcal{S}^i(0),\ {\rm dist}[x, \partial \mathcal{S}^i(0)] \geq \delta,\ i = 1,\dots, N.
		\end{array}\right. 		
	\label{r7}	
\end{align}	
	\end{itemize}

As shown in \cite{FeHiNe}, system \eqref{r1}--\eqref{r4}, supplemented 
with the initial conditions \eqref{r5} -- \eqref{r7}, admits 
a solution defined in $(0,T) \times \T^3$, see  
\cite{FeHiNe}. Moreover, performing successively the limits $\ep \to 0$, 
$\delta \to 0$, the solutions converge to the desired solution of the 
fluid structure interaction problem as long as there are no contacts of rigid objects, meaning up to the first contact. 

Our ultimate goal is therefore to show that the limit solution satisfies the 
energy balance in the differential form. To see this, we first record the energy balance for the original system:
\[
\frac{\D }{\dt} \int_{\T^3} \frac{1}{2} \vr |\vu|^2 \dx 
+ \int_{\T^3} [\mu]_\delta |\Ds \vu |^2 \dx + \int_{\T^3} \chi_\ep 
|\vu|^2 \dx \leq \int_{\T^3} \vr \vc{g} \cdot \vu \dx.
\]
As a matter of fact, due to the regularization of the convective terms, 
the approximate solutions are regular enough for the energy balance to be derived directly form the equations. Using ``distributional'' formulation, we may also include the initial energy, 
\begin{align} 
- \int_0^T \partial_t \psi \int_{\T^3} \frac{1}{2} \vr |\vu|^2 \dx \dt 
&+ \int_0^T \psi \int_{\T^3} [\mu]_\delta |\Ds \vu |^2 \dx \dt + \int_0^T \psi \int_{\T^3} \chi_\ep 
|\vu|^2 \dx \dt\br &\leq 
\psi(0) \int_{\T^3} \frac{1}{2} \vr_{0,\delta} |\vu_{0, \delta}|^2 \dx
+
\int_0^T \psi \int_{\T^3} \vr \vc{g} \cdot \vu \dx \dt
\label{r8}
\end{align}
for any $\psi \in C^1_c[0,T)$, $\psi \geq 0$.

Thanks to the specific choice of initial data, we have 
\begin{equation} \label{r9}
\int_{\T^3} \frac{1}{2} \vr_{0,\delta} |\vu_{0, \delta}|^2 \dx = 
\intO{ \frac{1}{2} \vr_{0,\delta} |\vu_{0, \delta}|^2 }.	
	\end{equation}
	
Next, because of the domain penalization term,
\begin{equation} \label{r10}
\int_0^T \int_{\T^3 \setminus \Omega} |\vu|^2 \aleq \ep\ 
\Rightarrow \ \int_0^T \int_{\T^3 \setminus \Omega} \vr |\vu|^2 \aleq \ep
	\end{equation}
since the densities remain bounded uniformly in time by their initial values.
Moreover, by the same token, the artificial viscosities $[\mu]_\delta$ remain bounded below, 
\begin{equation} \label{r11}
	[\mu]_{\delta} \geq 1.
	\end{equation}

In view of the above observations, we may rewrite \eqref{r8} as
\begin{align} 
	- \int_0^T \partial_t \psi \intO{ \frac{1}{2} \vr |\vu|^2 } \dt 
	&+ \int_0^T \psi \intO{ |\Ds \vu |^2 } \dt \br &\leq 
	\psi(0) \intO{ \frac{1}{2} \vr_{0,\delta} |\vu_{0, \delta}|^2 }
	+
	\int_0^T \psi \intO{ \vr \vc{g} \cdot \vu } \dt + 
	\omega\left(\ep, \| \psi \|_{C^1[0,T]} \right)
	\label{r12}
\end{align}
where the error term satisfies 
\[
\omega(\ep, \| \psi \|_{C^1[0,T]}) \to 0 \ \mbox{as}\ \ep \to 0,
\] 
for any (fixed) $\psi \in C^1_c[0,T)$, $\psi \geq 0$.

At  this stage, we can perform the double limit $\ep \to 0$, $\delta \to 0$ in \eqref{r12}. As shown in \cite{FeHiNe}, the convective term 
\begin{equation} \label{r13}
	\vr_{\ep, \delta} \vu_{\ep, \delta} \otimes 
	\vu_{\ep, \delta} \stackrel{\ep \to 0}{\to} \vr_\delta \vu_{\delta} \otimes \vu_{\delta} \stackrel{\delta \to 0}{\to} 
	\vr \vu \otimes \vu \ \mbox{weakly in}\ L^q((0,T) \times \Omega; R^{3 \times 3})\ \mbox{for some}\ q > 1.
\end{equation}
In particular, we conclude by the convergence of the weak time derivative of the kinetic energy,  
\[	 
\int_0^T \partial_t \psi \intO{ \frac{1}{2} \vr_{\ep, \delta} |\vu_{\ep, \delta}|^2 } \dt \stackrel{\ep \to 0}{\to} 
\int_0^T \partial_t \psi \intO{ \frac{1}{2} \vr_\delta |\vu_\delta|^2 }
\stackrel{\delta \to 0}{\to} \int_0^T \partial_t \psi \intO{ \frac{1}{2} \vr |\vu|^2 }.
\]

We have shown the following result.

\begin{Theorem}[{\bf Dissipative solution -- global existence}] \label{T1}	
	
Let $\Omega \subset R^3$ be a bounded domain. Let 
$(\mathcal{S}^i)_{i=1}^N$ be a family of connected compact sets satisfying \eqref{rr1}, \eqref{rr2}. 
Let the mass densities $(\vr_{\mathcal{S}^i})_{i=1}^N$ be positive constants, 
\[
0 < \underline{\vr} \leq \vr_{\mathcal{S}^i} \leq \Ov{\vr} 
\ \mbox{for all}\ i = 1,\dots, N.
\]	
Finally, let $\vc{g} \in L^\infty(\Omega; R^3)$ be a given volume force. 

Then the fluid--structure interaction problem admits a dissipative weak solution in $(0,\infty) \times \Omega$ in the sense specified in Definition \ref{D1} for any finite energy initial data.

\end{Theorem}


{\smallskip
\begin{Remark} \label{rmk:diss2d} 
One can obtain  global dissipative weak solutions in dimension two as well by using the construction proposed in \cite{SST} which provides global weak solutions satisfying an integral form of energy inequality. Their solutions are obtained in a similar manner as we do, by approximating with those of Navier-Stokes in the whole domain and  that are ``very viscous"  on the regions corresponding to the rigid bodies. These approximations can be checked to satisfy the energy balance in differential form \eqref{m13}. Then, based on the strong convergences obtained in \cite{SST}, namely $(3.13)-(3.15)$, one can pass to the limit in \eqref{m13} to obtain that the limit solution, the weak solution for the fluid-structure interaction system obtained in \cite{SST}, is also a dissipative weak solution in our sense. 
\end{Remark}
}

\section{Long--time behaviour of dissipative solutions}
\label{L}	

The fact the dissipative solutions satisfy the energy balance in the differential form \eqref{m13} has an immediate impact on their long--time
behaviour. 

\subsection{System without external forcing}

Consider the case $\vc{g} = 0$. First, by Korn--Poincar\' e inequality,
the energy balance \eqref{m13} gives rise to
\begin{align} 
	- \int_0^\infty \partial_t \psi \intO{ \frac{1}{2} \vr |\vu|^2 } \dt
	+ C_{K,P} \int_0^\infty \psi \intO{ |\vu|^2 } \dt \leq \frac{1}{2} \psi (0) \intO{ \vr_0 |\vu_0|^2 }
	\nonumber
\end{align}
for any $\psi \in C^1_c[0, \infty)$, $\psi \geq 0$ and a certain positive constant $C_{K,P}$. Moreover, as the density is uniformly bounded, we get 
\begin{equation} \label{L1}
- \int_0^\infty \partial_t \psi \intO{ \frac{1}{2} \vr |\vu|^2 } \dt
+ \frac{2 C_{K,P}}{\Ov{\vr}} \int_0^\infty \psi \intO{ \frac{1}{2} \vr |\vu|^2 } \dt \leq \frac{1}{2} \psi (0) \intO{ \vr_0 |\vu_0|^2 }.
\end{equation}
for any $\psi \in C^1_c[0, \infty)$, $\psi \geq 0$. Relation \eqref{L1} yields the exponential decay of the kinetic energy, 
\begin{equation} \label{L2}
\intO{ \vr |\vu|^2 (\tau, \cdot) } \leq \exp\left( - \frac{2 C_{K,P}}{\Ov{\vr}} \tau \right) \intO{ \vr_0 |\vu_0|^2 },\ \tau \geq 0. 
\end{equation}	

Finally, we use the fact the mass densities of the rigid objects are bounded below away from zero and their interiors are non--empty. This yields the exponential decay for the corresponding rigid velocities, 
\begin{equation} \label{L3}
| \vc{Y}_i(\tau) | = \left| \frac{\D }{\dt} \vc{h}_i (\tau) \right| \leq 
C \exp (- \Lambda \tau), \ |\mathbb{Q}_i(\tau)| \leq 
C \exp (- \Lambda \tau),\ i = 1, \dots, N, 
	\end{equation}
for some $\Lambda > 0$.

We have shown the following result. 

	
\begin{Theorem}[\bf Long--time behaviour without forcing] \label{LT1}
	
Under the hypotheses of Theorem \ref{T1}, suppose 
$\vc{g} = 0$. Then the kinetic energy of the system tends exponentially to zero as stated in \eqref{L2}. In particular, 
\[
\vc{h}_i(\tau, \cdot) \to \vc{h}_i^\infty, \ |\mathbb{Q}_i(\tau)| \to 0 
\ \mbox{as}\ \tau \to \infty, 
\]
for some points $\vc{h}^\infty_i$, $i = 1,\dots, N$.
	
\end{Theorem}	
	

{
\begin{Remark}
{Obviously, the dissipative solutions may not be uniquely determined by the initial data. We point out that Theorem \ref{LT1} applies to \emph{any}
global in time dissipative solution, not just to those constructed in the proof of Theorem \ref{T1}.}	 
Taking into account the previous Remark~\ref{rmk:diss2d} the result of Theorem~\ref{LT1} also applies to the two dimensional weak solutions obtained in \cite{SST}.

\end{Remark}
}
    
\subsection{Gravitational forcing} 

We conclude the paper by a short discussion of the asymptotic behaviour of the system driven by gravitational (potential) external force. 
Accordingly, we suppose 
\[
\vc{g} = \Grad G,\ G = G(x) \ \mbox{a smooth potential.}
\]
In this case, the forcing term on the right--hand side of \eqref{m13} can be written with the help of the equation of continuity as 
\[
\int_0^\infty \psi \intO{ \vr \Grad G \cdot \vu } = 
- \int_0^\infty \partial_t \psi \intO{ \vr G } + \psi(0) \intO{ \vr_0 G }.
\]	
Consequently, the energy inequality \eqref{m13} takes the form
\begin{align} 
	- \int_0^\infty \partial_t \psi \intO{ \vr \left( \frac{1}{2} |\vu|^2 - G \right) } \dt
	&+ \int_0^\infty \psi \intO{ \mathbb{S}(\Ds \vu) : \Ds \vu } \dt \br  
	&\leq \psi (0) \intO{ \vr_0 \left( \frac{1}{2} |\vu_0|^2 - G \right) }
	\label{L4}
\end{align} 	
for any $\psi \in C^1_c[0, \infty)$, $\psi \geq 0$.

The function 
\[
\tau \mapsto \intO{ \vr \left( \frac{1}{2} |\vu|^2 - G \right)(\tau, \cdot) }
\]
is non--increasing and bounded from below. We may infer 
\begin{equation} \label{L5}
\intO{ \vr \left( \frac{1}{2} |\vu|^2 - G \right)(\tau, \cdot) } \to 
- \mathcal{E}_\infty \ \mbox{as}\ \tau \to \infty.
\end{equation}

In addition, we have 
\begin{equation} \label{L6}
\int_0^\infty \intO{ |\Ds \vu |^2 } \dt < \infty.
\end{equation}
Moreover, by virtue of the equation of continuity , we deduce 
\begin{equation} \label{L7}
\intO{ \vr {G} (\tau_1, \cdot) } - 
\intO{ \vr {G} (\tau_2, \cdot) } = \int_{\tau_1}^{\tau_2} \intO{ 
\vr \vu \cdot \Grad G } \dt,
	\end{equation}
where, furthermore, 
\begin{align}
\left| \int_{\tau_1}^{\tau_2} \intO{ 
	\vr \vu \cdot \Grad G } \dt \right| &\leq 
\int_{\tau_1}^{\tau_2} \intO{ \vr |\Grad G|^2 } \dt +	
\int_{\tau_1}^{\tau_2} \intO{ \vr |\vu|^2 } \dt \br &\leq 
\int_{\tau_1}^{\tau_2} \intO{ \vr |\Grad G|^2 } \dt
+ C \int_{\tau_1}^{\tau_2} \intO{ |\Ds \vu |^2 } \dt.
\label{L8}
\end{align}

Combining \eqref{L5}--\eqref{L8} we conclude
\begin{equation} \label{L9}
\intO{ \vr |\vu|^2 (\tau, \cdot) } \to 0,\ 
\intO{ \vr(\tau, \cdot) G } \to \mathcal{E}_\infty \ \mbox{as}\ \tau \to \infty.	
	\end{equation}
We have shown the following result. 

	
	\begin{Theorem}[\bf Long--time behaviour, gravitational forcing] \label{LT2}
Under the hypotheses of Theorem \ref{T1}, suppose 
\[ 
\vc{g} = \Grad G,\ G \in C^1(\Ov{\Omega}).
\]

Then 
\[ 
\intO{ \vr |\vu|^2 (\tau, \cdot) } \to 0 \ \mbox{as}\ \tau \to \infty.
\]
Moreover, there is a constant $\mathcal{E}_\infty$ such that 	
\[
\intO{ \vr (\tau, \cdot) G } \to \mathcal{E}_\infty 
\ \mbox{as}\ \tau \to \infty.
\]
\end{Theorem}


In the full generality allowed by Theorem \ref{LT2}, it seems difficult 
to identify the final position of rigid bodies. More specific results in the simplified 2D case have been obtained in \cite{FeiNec2011}.
\bibliography{citace}
\bibliographystyle{plain}

\end{document}